\RequirePackage{xcolor}
\documentclass[journal,twoside,web]{ieeecolor}
\usepackage{generic}

\usepackage{amsmath,amssymb,amsfonts}
\usepackage{graphicx}
\usepackage{mathtools}
\usepackage{booktabs}
\usepackage{textcomp}

\usepackage[backend=bibtex,bibstyle=ieee,citestyle=numeric-comp,doi=false,isbn=false,eprint=true,url=false]{biblatex}
\newtheorem{theorem}{Theorem}
\newtheorem{proposition}{Proposition}
\newtheorem{lemma}{Lemma}
\newtheorem{remark}{Remark}
\newtheorem{assumption}{Assumption}

\usepackage{hyperref}
\hypersetup{
    colorlinks=true,
    linkcolor=blue,
    citecolor=blue,
    urlcolor=blue,
    bookmarksnumbered=true,
    bookmarksopen=true
}

\newcommand{\tr}{\operatorname{tr}}
\newcommand{\vect}{\operatorname{vec}}
\newcommand{\diag}{\operatorname{diag}}
\newcommand{\cl}{\operatorname{cl}}
\newcommand{\Sym}{\mathbb{S}}
\newcommand{\Nb}{\bar N}
\newcommand{\Kopt}{K_{\mathrm{opt}}}
\newcommand{\Kloc}{K_{\mathrm{loc}}}
\newenvironment{tbmatrix}{\setlength{\arraycolsep}{2.5pt}\begin{bmatrix}}{\end{bmatrix}}

\begin{document}

\title{A Lifted Hessian--Riemannian Gradient Flow for \\ LQR Synthesis under Passivity Constraints}

\author{Armin Gießler and Sören Hohmann, \IEEEmembership{Senior Member, IEEE}
\thanks{A. Gießler and S. Hohmann are with the Institute of Control Systems, Karlsruhe Institute of Technology (KIT), 76131 Karlsruhe, Germany. Corresponding author is Armin Gießler (e-mail: armin.giessler@kit.edu).}}

\maketitle

\begin{abstract}
The static state-feedback gain that renders a linear system passive at an external port and minimizes a linear quadratic regulator (LQR) cost cannot be computed by convex optimization, since the passivity certificate and the LQR cost are coupled bilinearly through the common gain. A previous indirect approach restricts the gain to a convex inner approximation of the passivity region, which is conservative and whose construction scales exponentially with the number of gain entries. This paper proposes a direct approach. The set of strictly passivating gains is the image of a convex lifted set of inverse storage matrices and transformed gains. We describe this set by a strict linear matrix inequality and an affine equality for every feedthrough matrix with positive semidefinite symmetric part, including feedthrough matrices with a singular symmetric part, and minimize the LQR cost over it by a gradient flow whose metric is the Hessian of a logarithmic barrier. For multi-input systems, this Hessian is singular along gain changes that preserve the storage function, and we complete the metric along these directions by the metric of the gain space. We prove that the flow exists for all times, that every gain along the flow is strictly passivating with an explicit storage function, and that the cost decreases monotonically while the gains remain in a compact set of stabilizing gains. Numerical examples show that the flow attains a lower cost than the indirect approach, that the completion is required for multi-input systems, and that the flow can converge to a local minimizer on nonconvex passivity regions.
\end{abstract}


\section{Introduction}
\label{sec:intro}

\IEEEPARstart{P}{assivity} describes the input-output behavior of a dynamical system by a dissipation inequality at its port \cite{willems1972dissipativeI,willems1972dissipativeII}. Feedback interconnections of passive systems inherit stability properties \cites[Section~6.5]{khalil_nonlinear_2002}[Sections~2.2 and 4.3]{van_der_schaft2017}[Chapter~2]{arcak2016networks}, where the precise statements depend on the type of interconnection, on the strictness of the passivity properties, and on detectability conditions. Passivity at an external port therefore permits stable interconnections with passive subsystems that are unknown a priori, but it does not ensure a sufficient transient performance. This motivates a state feedback that renders the closed-loop system passive at an external port and, among all such gains, minimizes a linear quadratic regulator (LQR) cost of the isolated system \cite{giessler2026passive}.

Each of the two requirements alone leads to a convex synthesis problem. For linear systems, the Kalman--Yakubovich--Popov (KYP) lemma characterizes passivity by a linear matrix inequality (LMI) in a storage matrix \cites[Section~2.7.2]{LMIbook}[Theorem~5]{Madeira2016}. For the closed loop, this inequality is bilinear in the storage matrix and the feedback gain, and a change of variables turns the synthesis of a feedback gain with a passivity requirement into an LMI problem \cite[Proposition~3]{giessler2026passive}. Likewise, the LQR problem admits convex reformulations through the covariance formulation \cite[Section~2]{Feron} or the Riccati inequality \cite[p.~115]{LMIbook}. If both requirements, passivity and LQR performance, are imposed jointly, however, they are coupled through the common gain, which results in a bilinear equality. Hence, the resulting formulations admit no convex reformulation \cite[Propositions~5 and~6]{giessler2026passive}, and the set of passivating gains is nonconvex in general \cite[Lemma~10]{giessler2026passive}. The indirect approach of \cite{giessler2026passive} therefore certifies cubes in the gain space by a common storage function, inner-approximates their union by a convex polytope, and minimizes the LQR cost over this polytope by a projected gradient flow. Its gain is suboptimal whenever the constrained minimizer lies outside the polytope, the number of cubes grows exponentially with the number of gain entries, and computing the largest inner polytope is open for more than two gain entries \cite[Remark~4]{giessler2026passive}.

This paper proposes a direct approach. The set of passivating gains is the image of a convex lifted set under the map from the inverse storage matrix and a transformed gain to the gain. The LQR cost is nonconvex in the lifted variables, but the lifted set admits a logarithmic barrier, which we use as the metric of a Hessian--Riemannian gradient flow \cite{alvarez2004}. The flow evolves the gain and its storage function jointly and generates passivating gains at all times.

\subsection{Literature Review}

\subsubsection*{Passivity-Based Feedback Design}
Dynamic output-feedback controllers that are required to be (strictly) positive real, which is equivalent to (strict) passivity for minimal realizations \cite[Section~II-B and Theorem~5]{Madeira2016}, are designed for LQG-type performance \cite{LozanoLeal1988,Jacobus1990}, for mixed $\mathcal{H}_2/\mathcal{H}_\infty$ performance \cite{Haddad1994}, and for the closed-loop $\mathcal{H}_2$-norm or an upper bound on it \cite{Geromel1997,shimomura2002strictly,Forbes2019,Forbes2013,Damaren2006}. Passivation of the closed loop without a performance criterion is addressed by static state feedback in \cite{byrnes_passivity_1991}, by robust static output feedback in \cite{peaucelle2005robust}, and by fixed-structure controllers in \cite{Su2020}. Necessary and sufficient conditions for dynamic output-feedback and static state-feedback controllers that render the closed loop extended strictly positive real are given in \cite{Sun1994}. Static feedback and an LQR-type criterion are combined in \cite{Hallinan}, where the gain is inverse optimal and the state-penalty matrix follows from the storage function, and in the indirect approach \cite{giessler2026passive}.

\subsubsection*{Policy Optimization}
Policy optimization minimizes the cost directly over the controller parameters \cite{Hu2023,talebi2024policyoptimizationcontrolgeometry}. For the LQR problem, the cost is coercive on the stabilizing gains and gradient dominated on its sublevel sets, so that gradient methods converge to the optimal gain \cite{fazel2018global,bu2020clqr,Mohammadi2019,Mohammadi2022,Bu2020}, and so does the gradient flow of a Bellman error parameterized in the gain \cite{giessler2025cdc}. 
Gradient dominance follows from convex parameterizations based on LMIs for a broad class of problems \cite{Fazel2021}, and convex liftings establish the global optimality of nondegenerate stationary points for LQG and $\mathcal{H}_\infty$ output feedback \cite{zheng_benign_2026}. Constraints on the policy are handled by Riemannian Newton-type methods on submanifolds of stabilizing controllers \cite{Talebi2022}, by projection \cite{Takakura2024}, or implicitly by the update \cite{Zhang2021}.

\subsubsection*{Riemannian Gradient Flows}
Riemannian gradient flows have been studied in connection with interior-point methods, see \cite{alvarez2004} and the references therein. For metrics given by the Hessian of a Legendre-type function, \cite{alvarez2004} establishes global existence and, for convex objective functions, convergence of the cost to its infimum. Logarithmic barriers of LMIs are the standard barriers of interior-point methods for semidefinite programming \cite{Nesterov1994}.

\subsection{Contributions}

\begin{enumerate}
    \item An explicit description of the lifted set of strictly passivating gains by a strict LMI and an affine equality, valid for feedthrough matrices with singular symmetric part (Theorem~\ref{thm:lifted}), and a first-order necessary condition in terms of the certificate (Proposition~\ref{prop:necessary}).
    \item A lifted Hessian--Riemannian gradient flow, completed for multi-input systems, with guaranteed global existence, strict passivity along the trajectory, and monotone cost decrease (Theorem~\ref{thm:flow}).
    \item A numerical comparison with the indirect approach \cite{giessler2026passive} and an analysis of the local convergence on nonconvex passivity regions.
\end{enumerate}

\subsection{Notation}
The sets of real $n \times m$ matrices and of symmetric $n \times n$ matrices are $\mathbb{R}^{n\times m}$ and $\Sym^n$. For $X \in \Sym^n$, $X \succ 0$ ($X \succeq 0$) denotes positive definiteness (semidefiniteness), and $\lambda_{\min}(X)$ and $\lambda_{\max}(X)$ denote the extreme eigenvalues. For a matrix $A$, $A^\top$, $\tr(A)$, $\vect(A)$, $\|A\|_F$, and $\|A\|_2$ denote the transpose, trace, columnwise vectorization, Frobenius norm, and spectral norm, and $\otimes$ is the Kronecker product. The identity matrix is $I_n$, its $i$-th column is $e_i$, and $\diag(\cdot)$ forms a (block) diagonal matrix. The interior, boundary, and closure of a set $\mathcal{X}$ are $\operatorname{int}(\mathcal{X})$, $\partial\mathcal{X}$, and $\cl(\mathcal{X})$. For a vector $v$, $\|v\|$ is the Euclidean norm and, for $G \succ 0$, $\|v\|_G = (v^\top G v)^{1/2}$.

\section{Preliminaries}
\label{sec:prelim}

\subsection{Passivity}
\label{subsec:passivity}

The LTI system
\begin{align}
    \dot x = A x + B_v v, \qquad y = C x + D v, \label{eq:sys_pass}
\end{align}
with $x \in \mathbb{R}^n$ and $v, y \in \mathbb{R}^p$ is strictly (state) passive if there exist a continuously differentiable, positive semidefinite storage function $S(x)$ and a positive definite function $\psi(x)$ such that
\begin{align}
    \dot S = \tfrac{\partial S}{\partial x}(x) \, (Ax + B_v v) \leq v^\top y - \psi(x) \label{eq:passivity}
\end{align}
for all $x \in \mathbb{R}^n$ and $v \in \mathbb{R}^p$, and it is passive if \eqref{eq:passivity} holds with $\psi \equiv 0$ \cite[Definition~6.3]{khalil_nonlinear_2002}. Passivity requires $D + D^\top \succeq 0$ \cite[Section~II-A and Appendix~A]{Madeira2016}. Define
\begin{align}
    M &\coloneqq \begin{bmatrix} A^\top P + PA & PB_v - C^\top \\ B_v^\top P - C & -D - D^\top \end{bmatrix}, \label{eq:M} \\
    N \succeq_s 0 \; &:\Leftrightarrow \; \begin{cases} N \succ 0, & \text{for strict passivity}, \\ N = 0, & \text{for nonstrict passivity}, \end{cases} \label{eq:succeq_s} \\
    L(D) &\coloneqq \begin{cases} A^\top P + PA + N \preceq 0, ~ B_v^\top P = C, & D = 0, \\ M + \diag(N,0) \preceq 0, & D \neq 0, \end{cases} \label{eq:L}
\end{align}
where $N \succeq_s 0$ in both cases, and let $L_{\mathrm{s}}(D)$ and $L_{\mathrm{ns}}(D)$ denote $L(D)$ with $N \succ 0$ and $N = 0$, respectively.

\begin{lemma}[KYP Lemma]
    \label{lem:KYP}
    If $L(D)$ holds for some $P \succ 0$, then \eqref{eq:sys_pass} is (strictly) passive with the storage function $S(x) = \tfrac{1}{2} x^\top P x$. If $(A,B_v)$ is controllable and $(C,A)$ is observable, the converse holds as well.
\end{lemma}

\begin{proof}
    The first statement follows for $D \neq 0$ from
\begin{align}
    \dot S &= v^\top y - \tfrac{1}{2} x^\top N x + \tfrac{1}{2} \begin{bmatrix} x \\ v \end{bmatrix}^\top \bigl(M + \diag(N,0)\bigr) \begin{bmatrix} x \\ v \end{bmatrix} \notag \\
    &\leq v^\top y - \tfrac{1}{2} x^\top N x
\end{align}
and analogously for $D = 0$. The converse is the positive real lemma \cites[Lemma~2]{kottenstette2014relationships}[Section~2.7.2]{LMIbook} for the nonstrict case and follows from \cite[Section~III and Theorem~5]{Madeira2016} for the strict case, where, for $D = 0$, the zero lower-right block of $M$ forces $B_v^\top P = C$ \cite[Theorem~1]{Albert1969}.
\end{proof}

Since $N$ enters only the state block, Lemma~\ref{lem:KYP} admits singular $D + D^\top$, whereas the strict inequality $M \prec 0$ used for $D \neq 0$ in \cite{giessler2026passive} requires $D + D^\top \succ 0$ \cite[Lemma~2.3]{Sun1994}.

\subsection{LQR Cost and Its Gradient}
\label{subsec:lqr}

For $\dot x = Ax + B_u u$ with $u = -Kx$, let $A_K = A - B_u K$ and $\mathcal{K} = \{ K \in \mathbb{R}^{m \times n} \mid A_K \text{ is Hurwitz} \}$. For $Q \succeq 0$ and $R \succ 0$, the LQR cost is $f_K\colon \mathcal{K} \to \mathbb{R}_{\geq 0}$, $K \mapsto \tr(X_K)$, where $X_K \succeq 0$ is the unique solution of
\begin{align}
   A_K^\top X_K + X_K A_K + Q + K^\top R K = 0. \label{eq:cost}
\end{align}
It is the sum of the costs $\int_0^\infty x^\top Q x + u^\top R u \, \mathrm{d}t$ over the initial states $e_1, \dots, e_n$. If $(A,B_u)$ is stabilizable and $(Q^{1/2},A)$ is detectable, then this sum is infinite for every $K \notin \mathcal{K}$ \cite[Section~3.3]{bu2020clqr}, $f_K$ is real analytic on $\mathcal{K}$, $f_K \to \infty$ as $K \to \partial\mathcal{K}$ or $\|K\|_F \to \infty$, and the sublevel sets $\{K \in \mathcal{K} \mid f_K \leq c\}$ are compact \cite[Section~3.4, Lemma~3.3, and Corollary~3.3.1]{bu2020clqr}. The gradient of $f_K$ is \cites[Section~IV]{Levine1970}[Proposition~3.4]{bu2020clqr}
\begin{align}
    \nabla f_K = 2 (R K - B_u^\top X_K) Y_K, \label{eq:gradK}
\end{align}
where $Y_K \succ 0$ is the unique solution of $A_K Y_K + Y_K A_K^\top + I_n = 0$, and the optimal LQR gain $K^*$ is the unique stationary point of $f_K$ on $\mathcal{K}$ \cite[Lemma~3.5 and Corollary~3.5.1]{bu2020clqr}.

\subsection{Hessian--Riemannian Gradient Flows}
\label{subsec:hr}

Let $\Theta \subseteq \mathbb{R}^{n_\theta}$ be open and convex, let $F \in C^1(\Theta)$, and let $h \in C^2(\Theta)$ with $G(\theta) = \nabla^2 h(\theta) \succ 0$ for all $\theta \in \Theta$. The gradient of $F$ with respect to the metric $(u,v)_\theta = u^\top G(\theta) v$ is $G(\theta)^{-1} \nabla F(\theta)$, and
\begin{align}
    \dot\theta = -G(\theta)^{-1} \nabla F(\theta) \label{eq:hr}
\end{align}
is referred to as a Hessian--Riemannian gradient flow \cite{alvarez2004}. Along \eqref{eq:hr}, the chain rule gives $\tfrac{\mathrm{d}}{\mathrm{d}t} F(\theta) = -\|\nabla F(\theta)\|^2_{G(\theta)^{-1}}$. If $h$ is a barrier of $\Theta$, namely $h(\theta) \to \infty$ as $\theta \to \partial\Theta$, such as the logarithmic barrier $h(\theta) = -\sum_i \log\theta_i$ of $\Theta = (0,\infty)^{n_\theta}$, then $h$ is of Legendre type \cite[Definition~3.1 and Section~4.4]{alvarez2004}. For such $h$ and $F$ continuously differentiable on $\mathbb{R}^{n_\theta}$ with a locally Lipschitz continuous gradient, global existence of the solution of \eqref{eq:hr} in $\Theta$ and, if $F$ is also convex, convergence of $F(\theta(t))$ to $\inf_\Theta F$ are established under additional conditions in \cite[Theorem~4.1 and Proposition~4.1]{alvarez2004}.

\section{Problem Formulation}
\label{sec:problem}

Consider the LTI system
\begin{align}
    \dot x = A x + B_u u + B_d d, \qquad y = C x + D d, \label{eq:sys}
\end{align}
with $x \in \mathbb{R}^n$, $u \in \mathbb{R}^m$, and $d, y \in \mathbb{R}^p$. The state feedback $u = -Kx$ yields the closed-loop system
\begin{align}
    \dot x = A_K x + B_d d, \qquad y = C x + D d, \label{eq:closed}
\end{align}
with the external port $(d,y)$, at which passivity is imposed. Let $M_K$ and $L(D,K)$ denote $M$ and $L(D)$ from \eqref{eq:M} and \eqref{eq:L} with $A$ and $B_v$ replaced by $A_K$ and $B_d$, respectively. 
The sets of gains that render \eqref{eq:closed} passive and strictly passive are
\begin{align}
    \mathcal{P}^0 &\coloneqq \{ K \in \mathbb{R}^{m \times n} \mid \exists P \succ 0 \colon L_{\mathrm{ns}}(D,K) \}, \label{eq:P0} \\
    \mathcal{P}^+ &\coloneqq \{ K \in \mathbb{R}^{m \times n} \mid \exists P \succ 0, N \succ 0 \colon L_{\mathrm{s}}(D,K) \}, \label{eq:Pplus}
\end{align}
and are referred to as the passivity region and the strict passivity region, respectively. The symbol $\mathcal{P}^{(+)} \in \{\mathcal{P}^0, \mathcal{P}^+\}$ refers to either of them, and a matrix $P \succ 0$ for which $L(D,K)$ holds is called a certificate of $K$. By Lemma~\ref{lem:KYP}, every $K \in \mathcal{P}^{(+)}$ renders \eqref{eq:closed} (strictly) passive, and the converse holds if $(A_K,B_d)$ is controllable and $(C,A_K)$ is observable.

\begin{assumption}
    \label{ass:sys}
    The pair $(A,B_u)$ is stabilizable, $(Q^{1/2},A)$ is detectable, $B_u$ has full column rank, and $\mathcal{P}^+ \neq \emptyset$.
\end{assumption}

The leading block of $L(D,K)$ gives $A_K^\top P + PA_K + N \preceq 0$, so that $\mathcal{P}^+ \subseteq \mathcal{K}$ and $\mathcal{P}^+ \subseteq \mathcal{P}^0 \subseteq \cl(\mathcal{K})$ \cite[Lemma~4]{giessler2026passive}, and for $D \neq 0$, its lower-right block gives $D + D^\top \succeq 0$, so that $\mathcal{P}^{(+)} = \emptyset$ otherwise. We seek a solution of
\begin{align}
    \min_{K} ~ f_K \quad \text{s.t.} \quad K \in \mathcal{K} \cap \mathcal{P}^{(+)}. \label{eq:problem}
\end{align}

\begin{proposition}
    \label{prop:feas}
    It holds that $K \in \mathcal{P}^{(+)}$ if and only if $K = W Z^{-1}$ for some $Z \succ 0$, $W \in \mathbb{R}^{m \times n}$, and $\tilde N \succeq_s 0$ such that, for $D \neq 0$,
    \begin{align}
        N_D(Z,W,\tilde N) \coloneqq \begin{bmatrix} N_0(Z,W,\tilde N) & B_d - ZC^\top \\ B_d^\top - CZ & -D - D^\top \end{bmatrix} \preceq 0, \label{eq:ND}
    \end{align}
    and, for $D = 0$,
    \begin{align}
        N_0(Z,W,\tilde N) &\coloneqq ZA^\top + AZ - W^\top B_u^\top - B_u W + \tilde N \preceq 0, \notag \\
        B_d^\top &= CZ. \label{eq:N0}
    \end{align}
    The storage function is then $S(x) = \tfrac{1}{2} x^\top Z^{-1} x$.
\end{proposition}

\begin{proof}
    With $Z = P^{-1}$, $W = KZ$, and $\tilde N = ZNZ$, the congruence transformation with $\diag(Z,I_p)$ turns $M_K + \diag(N,0) \preceq 0$ into \eqref{eq:ND}, and for $D = 0$, the congruence transformation with $Z$ and the right multiplication of $B_d^\top P = C$ by $Z$ yield \eqref{eq:N0}. All steps are reversible, and $\tilde N \succeq_s 0$ if and only if $N \succeq_s 0$.
\end{proof}

For $D = 0$ or $D + D^\top \succ 0$, Proposition~\ref{prop:feas} reduces to \cite[Proposition~3]{giessler2026passive}. Let $\mathcal{L}_D^{(+)}$ and $\mathcal{L}_0^{(+)}$ denote the sets of $(Z,W) \in \Sym^n \times \mathbb{R}^{m \times n}$ with $Z \succ 0$ for which \eqref{eq:ND} and \eqref{eq:N0}, respectively, hold for some $\tilde N \succeq_s 0$, let $\mathcal{L}^{(+)}$ denote $\mathcal{L}_D^{(+)}$ for $D \neq 0$ and $\mathcal{L}_0^{(+)}$ for $D = 0$, and let $\varphi(Z,W) \coloneqq W Z^{-1}$. By Proposition~\ref{prop:feas}, $\mathcal{P}^{(+)} = \varphi(\mathcal{L}^{(+)})$, where $\mathcal{L}^{(+)}$ is convex as the projection of a convex set.

The indirect approach \cite{giessler2026passive} certifies cubes $\mathcal{C}_i \subseteq \mathcal{P}^{(+)}$ by a common storage function at their vertices \cite[Lemma~11]{giessler2026passive}, inner-approximates their union by a convex polytope $\mathcal{P}_C$ \cite[Lemma~12]{giessler2026passive}, and computes $\hat K = \arg\min_{K \in \mathcal{K} \cap \mathcal{P}_C} f_K$ by a projected gradient flow \cite[Theorem~14]{giessler2026passive}. The gain $\hat K$ is feasible for \eqref{eq:problem}, but suboptimal whenever the minimizers of \eqref{eq:problem} lie outside the approximation $\mathcal{P}_C$.

\section{Lifted Formulation}
\label{sec:lifted}

\subsection{Explicit Description of the Lifted Set}

The logarithmic barrier, and thus the metric of \eqref{eq:hr}, is only defined if the interior of the feasible set is characterized by a strict LMI. Since $\tilde N$ enters only the state block, however, \eqref{eq:ND} is never strictly negative definite if $D + D^\top$ is singular, so that a barrier of $\mathcal{L}^+$ requires a description without $\tilde N$. For $D + D^\top \succeq 0$, let the columns of $U \in \mathbb{R}^{p \times r}$ and $V \in \mathbb{R}^{p \times (p - r)}$ be orthonormal bases of the range and the kernel of $D + D^\top$, where $r = \operatorname{rank}(D + D^\top)$, let $\Lambda \coloneqq U^\top (D + D^\top) U \succ 0$, and define
\begin{align}
   \!\! \Nb(Z,W) \! \coloneqq \begin{bmatrix} N_0(Z,W,0) & \!\!(B_d - ZC^\top) U \\ U^\top (B_d^\top - CZ)\!\! & -\Lambda \end{bmatrix}\! \in \Sym^{n+r}. \label{eq:Nbar}
\end{align}
\begin{theorem}
    \label{thm:lifted}
    Let $D + D^\top \succeq 0$. Then
    \begin{align}
        \mathcal{L}^+ &= \{ (Z,W) \mid Z \succ 0, ~ V^\top C Z = V^\top B_d^\top, ~ \Nb \prec 0 \}, \label{eq:Lplus} \\
        \mathcal{L}^0 &= \{ (Z,W) \mid Z \succ 0, ~ V^\top C Z = V^\top B_d^\top, ~ \Nb \preceq 0 \}, \label{eq:L0}
    \end{align}
    namely, $\mathcal{L}^+$ is described by a strict LMI in $(Z,W)$ on an affine subspace, without $\tilde N$. Moreover, $\mathcal{P}^+$ is open.
\end{theorem}

\begin{proof}
    For $D = 0$, $U$ is empty, $V = I_p$, and $\Nb(Z,W) = N_0(Z,W,0)$. A matrix $\tilde N \succ 0$ with $N_0(Z,W,0) + \tilde N \preceq 0$ exists if and only if $N_0(Z,W,0) \prec 0$, which shows \eqref{eq:Lplus}, and \eqref{eq:L0} holds with $\tilde N = 0$. For $D \neq 0$, let $\Omega = B_d - ZC^\top$. The congruence transformation with the orthogonal matrix $\diag(I_n, [U ~ V])$ turns \eqref{eq:ND} into
    \begin{align}
        \begin{bmatrix} N_0(Z,W,\tilde N) & \Omega U & \Omega V \\ U^\top \Omega^\top & -\Lambda & 0 \\ V^\top \Omega^\top & 0 & 0 \end{bmatrix} \preceq 0.
    \end{align}
    By \cite[Theorem~1]{Albert1969}, the zero diagonal block forces $\Omega V = 0$, which is $V^\top C Z = V^\top B_d^\top$, and the remaining inequality is equivalent to $N_0(Z,W,0) + \Omega U \Lambda^{-1} U^\top \Omega^\top \preceq -\tilde N$ by the Schur complement \cite[Section~A.5.5]{boyd2004convex}, since $\Lambda \succ 0$. A matrix $\tilde N \succ 0$ with this property exists if and only if the left-hand side is negative definite, which is equivalent to $\Nb(Z,W) \prec 0$ by the Schur complement, and $\tilde N = 0$ yields $\Nb(Z,W) \preceq 0$. Finally, if $(Z, KZ) \in \mathcal{L}^+$, then $(Z, K'Z) \in \mathcal{L}^+$ for all $K'$ close to $K$, since $W$ enters only $\Nb$, so that $\mathcal{P}^+$ is open.
\end{proof}

Theorem~\ref{thm:lifted} unifies three cases. For $D + D^\top = 0$, which includes $D = 0$, $\Nb = N_0(Z,W,0)$ is accompanied by the equality $CZ = B_d^\top$ in all $p$ output directions. For $D + D^\top \succ 0$, $V$ is empty and $\Nb$ is congruent to $N_D(Z,W,0)$, which recovers the strict condition of \cite[Proposition~3]{giessler2026passive}. For singular $D + D^\top \neq 0$, a reduced equality in the $p - r$ directions of the kernel and a bordered inequality of size $n + r < n + p$ remain. The set \eqref{eq:Lplus} does not depend on the choice of $U$ and $V$, since another orthonormal basis $UO$ with orthogonal $O$ transforms $\Nb$ by a congruence with $\diag(I_n, O)$, and another basis of the kernel leaves the equality unchanged.

\subsection{Strict and Nonstrict Passivity}
\label{subsec:strict}

\begin{lemma}
    \label{lem:closure}
    Let Assumption~\ref{ass:sys} hold. Then $\mathcal{P}^0 \subseteq \cl(\mathcal{P}^+)$ and
    \begin{align}
        \inf_{K \in \mathcal{P}^+} f_K = \inf_{K \in \mathcal{K} \cap \mathcal{P}^0} f_K. \label{eq:inf}
    \end{align}
    If $K^* \notin \mathcal{P}^+$, then $f_K$ has no local minimizer in $\mathcal{P}^+$, and every minimizer of $f_K$ over $\mathcal{K} \cap \mathcal{P}^0$ lies in $\partial \mathcal{P}^+$.
\end{lemma}

\begin{proof}
    Let $(Z', W') \in \mathcal{L}^0$ and $(Z'', W'') \in \mathcal{L}^+$. For $\mu \in (0,1]$, the convex combination $(Z_\mu, W_\mu) = (1-\mu)(Z',W') + \mu (Z'',W'')$ satisfies $Z_\mu \succ 0$ and $V^\top C Z_\mu = V^\top B_d^\top$, and, since $\Nb$ is affine, $\Nb(Z_\mu, W_\mu) = (1-\mu)\Nb(Z',W') + \mu\Nb(Z'',W'') \prec 0$. Hence, $(Z_\mu, W_\mu) \in \mathcal{L}^+$ by Theorem~\ref{thm:lifted}. Since $\varphi$ is continuous at $(Z', W')$, where $Z' \succ 0$, $\varphi(Z_\mu, W_\mu) \to \varphi(Z',W')$ as $\mu \to 0$, which shows $\mathcal{P}^0 \subseteq \cl(\mathcal{P}^+)$. Since $\mathcal{P}^+ \subseteq \mathcal{K} \cap \mathcal{P}^0$, the left-hand side of \eqref{eq:inf} is not smaller than the right-hand side. Conversely, every $K \in \mathcal{K} \cap \mathcal{P}^0$ is the limit of gains in $\mathcal{P}^+$, and $f_K$ is continuous on $\mathcal{K}$, which shows \eqref{eq:inf}. Since $\mathcal{P}^+$ is open, a local minimizer of $f_K$ in $\mathcal{P}^+$ is a stationary point and thus equal to $K^*$. A minimizer over $\mathcal{K} \cap \mathcal{P}^0$ lies in $\cl(\mathcal{P}^+)$, and if it were in $\mathcal{P}^+$, it would be a local minimizer in $\mathcal{P}^+$.
\end{proof}

By Lemma~\ref{lem:closure}, strictly passivating gains approximate the optimal value of the nonstrict problem arbitrarily well, and if $K^* \notin \mathcal{P}^+$, a method that generates strictly passivating gains can approach but not attain the minimizers of \eqref{eq:problem} for $\mathcal{P}^0$. If $K^* \in \mathcal{P}^+$, $K^*$ solves \eqref{eq:problem}.

\subsection{Certificate-Preserving Gain Changes}

A storage matrix certifies not only a single gain but a whole cone of gains.

\begin{lemma}
    \label{lem:cone}
    Let $L(D,K)$ hold with $P \succ 0$ and $N \succeq_s 0$, and let $\Sigma \in \mathbb{R}^{m \times m}$ with $\Sigma + \Sigma^\top \succeq 0$. Then $L(D, K + \Sigma B_u^\top P)$ holds with the same $P$ and $N$. If $\Sigma$ is skew-symmetric, then $M_{K + \Sigma B_u^\top P} = M_K$.
\end{lemma}

\begin{proof}
    For $K' = K + \Sigma B_u^\top P$, it holds that $A_{K'}^\top P + P A_{K'} = A_K^\top P + P A_K - P B_u (\Sigma + \Sigma^\top) B_u^\top P \preceq A_K^\top P + P A_K$, with equality for skew-symmetric $\Sigma$, whereas the remaining blocks of $M_K$ and the equality $B_d^\top P = C$ do not depend on $K$. Hence, $L(D,K')$ holds with the same $P$ and $N$.
\end{proof}

For $\Sigma = \kappa I_m$ with $\kappa \geq 0$, Lemma~\ref{lem:cone} yields the unbounded direction of $\mathcal{P}^{(+)}$ used in \cite[Lemma~7]{giessler2026passive}. For $m \geq 2$, the gain changes $\Sigma B_u^\top P$ with skew-symmetric $\Sigma$ form the subspace $\mathcal{T}(P) \coloneqq \{ \Sigma B_u^\top P \mid \Sigma^\top = -\Sigma \}$ of dimension $m(m-1)/2$, along which the gain changes without changing the certificate. This yields a first-order necessary condition for \eqref{eq:problem} that involves the certificate.

\begin{proposition}
    \label{prop:necessary}
    Let $\bar K$ be a local minimizer of $f_K$ over $\mathcal{K} \cap \mathcal{P}^0$, and let $L_{\mathrm{ns}}(D,\bar K)$ hold with $\bar P \succ 0$. Then $\nabla f_{\bar K} \bar P B_u$ is symmetric and positive semidefinite.
\end{proposition}

\begin{proof}
    Let $\Sigma + \Sigma^\top \succeq 0$ and $\Delta K \coloneqq \Sigma B_u^\top \bar P$. By Lemma~\ref{lem:cone} and since $\mathcal{K}$ is open, $\bar K + \epsilon \Delta K \in \mathcal{K} \cap \mathcal{P}^0$ for all sufficiently small $\epsilon > 0$. Since $\bar K$ is a local minimizer, the derivative of $f_K$ at $\bar K$ in the direction $\Delta K$ is nonnegative, namely
    \begin{align*}
        0 \leq \lim_{\epsilon \downarrow 0} \tfrac{f_{\bar K + \epsilon \Delta K} - f_{\bar K}}{\epsilon} = \tr(\nabla f_{\bar K}^\top \Delta K) = \tr(\Sigma^\top \nabla f_{\bar K} \bar P B_u).
    \end{align*}
    Skew-symmetric $\Sigma$ and $-\Sigma$ show that $\nabla f_{\bar K} \bar P B_u$ is symmetric, and symmetric $\Sigma \succeq 0$ show that it is positive semidefinite.
\end{proof}

For $m = 1$, Proposition~\ref{prop:necessary} reduces to $\nabla f_{\bar K} \bar P B_u \geq 0$. For $m \geq 2$, its symmetry part consists of $m(m-1)/2$ equations that a minimizer and its certificate $\bar P$ satisfy jointly.

\subsection{Coordinates and Lifted Objective}

The symmetric solutions of $V^\top C Z = V^\top B_d^\top$ form an affine subspace, which is nonempty by Assumption~\ref{ass:sys} and is parameterized by $Z(\xi) = Z_p + \sum_{j=1}^q \xi_j Z_j$ with a particular solution $Z_p$ and a basis $Z_1, \dots, Z_q$ of the homogeneous solutions.\footnote{The homogeneous solutions are the symmetric matrices whose range lies in the kernel of $V^\top C$, so that $q = (n - \operatorname{rank}(V^\top C))(n - \operatorname{rank}(V^\top C) + 1)/2$, where $\operatorname{rank}(V^\top C) = p - r$ if $C$ has full row rank.} With $\theta = (\xi, \vect(W))  \in \mathbb{R}^{n_\theta}$, $n_\theta\! = q + mn$, let
\begin{align}
    \Theta \coloneqq \{ \theta \in \mathbb{R}^{n_\theta} \mid (Z(\xi), W) \in \mathcal{L}^+ \}. \label{eq:Theta}
\end{align}

\begin{lemma}
    \label{lem:Theta}
    Let Assumption~\ref{ass:sys} hold. Then $\Theta = \{ \theta \in \mathbb{R}^{n_\theta} \mid Z(\xi) \succ 0, ~ \Nb(Z(\xi),W) \prec 0 \}$, and $\Theta$ is nonempty, open, and convex.
\end{lemma}

\begin{proof}
    Since $Z(\xi)$ satisfies the equality in \eqref{eq:Lplus} for all $\xi$, the description of $\Theta$ follows from Theorem~\ref{thm:lifted}, and $\Theta$ is open and convex, since $Z(\xi)$ and $\Nb(Z(\xi),W)$ are affine in $\theta$. By Assumption~\ref{ass:sys} and Proposition~\ref{prop:feas}, $\mathcal{L}^+ \neq \emptyset$, and every $(Z,W) \in \mathcal{L}^+$ satisfies $Z = Z(\xi)$ for some $\xi$ by \eqref{eq:Lplus}, so that $\Theta \neq \emptyset$.
\end{proof}

The partial derivatives $\tfrac{\partial Z}{\partial \theta_i}$ and $\tfrac{\partial \Nb}{\partial \theta_i}$ are constant, i.e.,
\begin{align}
    \frac{\partial \Nb}{\partial \xi_j} &= \begin{bmatrix} Z_j A^\top + A Z_j & -Z_j C^\top U \\ -U^\top C Z_j & 0 \end{bmatrix}, \notag \\
    \frac{\partial \Nb}{\partial W_{ab}} &= \begin{bmatrix} -B_u E_{ab} - E_{ab}^\top B_u^\top & 0 \\ 0 & 0 \end{bmatrix}, \label{eq:derivs}
\end{align}
$\partial Z / \partial \xi_j = Z_j$, and $\partial Z / \partial W_{ab} = 0$, where $E_{ab} \in \mathbb{R}^{m \times n}$ has a single one in entry $(a,b)$. The lifted objective is
\begin{align}
    F(\theta) \coloneqq f_{\varphi(Z(\xi),W)}, \qquad \theta \in \Theta, \label{eq:F}
\end{align}
which is well defined, since $\varphi(Z(\xi),W) \in \mathcal{P}^+ \subseteq \mathcal{K}$.

\begin{lemma}
    \label{lem:gradF}
    The function $F$ is real analytic and nonnegative on $\Theta$. With $K = W Z^{-1}$ and $Z = Z(\xi)$, its gradient is
    \begin{align}
        \nabla_W F = \nabla f_K \, Z^{-1}, \qquad \frac{\partial F}{\partial \xi_j} = -\tr\bigl(Z^{-1} \nabla f_K^\top K Z_j\bigr), \label{eq:gradF}
    \end{align}
    and $\nabla F(\theta) = 0$ if and only if $K = K^*$.
\end{lemma}

\begin{proof}
    $F$ is the composition of $f_K$, which is real analytic on $\mathcal{K}$, and of the real analytic map $\theta \mapsto W Z(\xi)^{-1}$. The differential $\mathrm{d}K = \mathrm{d}W \, Z^{-1} - K \, \mathrm{d}Z \, Z^{-1}$, which follows from $\mathrm{d}(Z^{-1}) = -Z^{-1}\,\mathrm{d}Z\,Z^{-1}$, with $\mathrm{d}Z = \sum_j Z_j\,\mathrm{d}\xi_j$ and $\mathrm{d}f_K = \tr(\nabla f_K^\top \mathrm{d}K)$ yield \eqref{eq:gradF}. Since $Z$ is invertible, $\nabla_W F = 0$ if and only if $\nabla f_K = 0$, namely $K = K^*$, which also implies $\partial F / \partial \xi_j = 0$.
\end{proof}

The function $F$ is nonconvex in general, and its Euclidean gradient $\nabla_W F$ scales with $Z^{-1}$, so that it becomes small for large $Z$ although $\nabla f_K$ does not.

\section{Lifted Hessian--Riemannian Gradient Flow}
\label{sec:flow}

\subsection{Barrier Metric}

The barrier
\begin{align}
    \mathcal{B}(\theta) \coloneqq -\log\det\bigl(-\Nb(Z(\xi),W)\bigr) - \log\det Z(\xi) \label{eq:barrier}
\end{align}
is convex on $\Theta$ and tends to infinity as $\theta$ approaches any point of $\partial\Theta$, at which $Z(\xi)$ or $\Nb(Z(\xi),W)$ is singular. Since $\Nb$ and $Z$ are affine in $\theta$, the identities $\mathrm{d}\log\det X = \tr(X^{-1}\,\mathrm{d}X)$ and $\mathrm{d}(X^{-1}) = -X^{-1}\,\mathrm{d}X\,X^{-1}$ for $X \succ 0$ show that, with $\Gamma \coloneqq (-\Nb)^{-1}$, the Hessian $H(\theta)$ of $\mathcal{B}$ has the entries
\begin{align}
    H_{ik}(\theta) = \tr\bigl(\Gamma \tfrac{\partial \Nb}{\partial \theta_i} \Gamma \tfrac{\partial \Nb}{\partial \theta_k}\bigr) + \tr\bigl(Z^{-1} \tfrac{\partial Z}{\partial \theta_i} Z^{-1} \tfrac{\partial Z}{\partial \theta_k}\bigr), \label{eq:hessB}
\end{align}
or, with $J_N = [\vect(\partial \Nb / \partial \theta_1), \dots, \vect(\partial \Nb / \partial \theta_{n_\theta})]$ and $J_Z = [\vect(\partial Z / \partial \theta_1), \dots, \vect(\partial Z / \partial \theta_{n_\theta})]$,
\begin{align}
    H(\theta) = J_N^\top (\Gamma \otimes \Gamma) J_N + J_Z^\top (Z^{-1} \otimes Z^{-1}) J_Z \succeq 0. \label{eq:hessvec}
\end{align}

\begin{lemma}
    \label{lem:kernel}
    Let Assumption~\ref{ass:sys} hold. For all $\theta \in \Theta$, the kernel of $H(\theta)$ is
    \begin{align}
        \mathcal{S} \coloneqq \bigl\{ (0, \vect(\Sigma B_u^\top)) \bigm| \Sigma \in \mathbb{R}^{m \times m}, ~ \Sigma^\top = -\Sigma \bigr\}, \label{eq:S}
    \end{align}
    which has dimension $m(m-1)/2$. In particular, $H(\theta) \succ 0$ if and only if $m = 1$.
\end{lemma}

\begin{proof}
    By \eqref{eq:hessvec}, $v = (v_\xi, \vect(\Delta W))$ with $v_\xi \in \mathbb{R}^q$ and $\Delta W \in \mathbb{R}^{m \times n}$ lies in the kernel of $H(\theta)$ if and only if $J_Z v = 0$ and $J_N v = 0$. The first condition is $v_\xi = 0$, since $Z_1, \dots, Z_q$ are linearly independent, and then the second condition is $B_u \Delta W + \Delta W^\top B_u^\top = 0$ by \eqref{eq:derivs}, namely that $B_u \Delta W$ is skew-symmetric. In this case, the range of $B_u \Delta W$ and the range of its transpose lie in the range of $B_u$, so that $B_u \Delta W = B_u \Sigma B_u^\top$ with the skew-symmetric matrix $\Sigma = (B_u^\top B_u)^{-1} B_u^\top (B_u \Delta W) B_u (B_u^\top B_u)^{-1}$, and $\Delta W = \Sigma B_u^\top$ since $B_u$ has full column rank. Conversely, $B_u \Sigma B_u^\top$ is skew-symmetric for every skew-symmetric $\Sigma$, and $\Sigma \mapsto \Sigma B_u^\top$ is injective.
\end{proof}

For $m \geq 2$, the Hessian $H$ is singular and thus does not define a metric for \eqref{eq:hr}. Along a direction $v = (0, \vect(\Sigma B_u^\top)) \in \mathcal{S}$, neither $Z$ nor $\Nb$ changes, so that the barrier is constant and the line $\theta + \alpha v$, $\alpha \in \mathbb{R}$, lies in $\Theta$. The gain, however, changes by $\alpha\Sigma B_u^\top Z^{-1} \in \mathcal{T}(Z^{-1})$, which keeps the certificate by Lemma~\ref{lem:cone} but changes $f_K$, and by Proposition~\ref{prop:necessary}, a minimizer is stationary along $\mathcal{T}$. We therefore complete $H$ by measuring these directions by the change of the gain. Let $J_K(\theta)$ denote the Jacobian of $\theta \mapsto \vect(W Z(\xi)^{-1})$, let the projection matrix $\Pi_{\mathcal{T}}(\theta) \in \Sym^{mn}$ denote the orthogonal projection onto $\vect(\mathcal{T}(Z(\xi)^{-1}))$, and define
\begin{align}
    G(\theta) \coloneqq H(\theta) + J_K(\theta)^\top \Pi_{\mathcal{T}}(\theta) J_K(\theta). \label{eq:G}
\end{align}

\begin{lemma}
    \label{lem:metric}
    Let Assumption~\ref{ass:sys} hold. Then, for all $\theta \in \Theta$,
    \begin{enumerate}
        \item[(i)] $G(\theta) \succeq H(\theta)$ and $G(\theta) \succ 0$, where $G = H$ for $m = 1$,
        \item[(ii)] $\|\nabla F(\theta)\|_{G(\theta)^{-1}} \geq \|\Pi_{\mathcal{T}}(\theta) \vect(\nabla f_K)\|$ with $K = W Z(\xi)^{-1}$,
        \item[(iii)] under every affine reparameterization $\theta = T \tilde\theta + c$ with invertible $T$, the metric in the coordinates $\tilde\theta$ is $T^\top G(T\tilde\theta + c) T$.
    \end{enumerate}
\end{lemma}

\begin{proof}
    (i) It holds that $v^\top G v = v^\top H v + \|\Pi_{\mathcal{T}} J_K v\|^2$. For $v \notin \mathcal{S}$, the first term is positive by Lemma~\ref{lem:kernel}. For $v = (0, \vect(\Sigma B_u^\top)) \in \mathcal{S} \setminus \{0\}$, $J_K v = \vect(\Sigma B_u^\top Z^{-1})$ is a nonzero element of $\vect(\mathcal{T})$, so that the second term is positive. For $m = 1$, $\mathcal{T} = \{0\}$. (ii) Since $J_K$ maps $\mathcal{S}$ onto $\vect(\mathcal{T})$, there exists $v \in \mathcal{S}$ with $J_K v = \Pi_{\mathcal{T}} \vect(\nabla f_K)$, for which $\nabla F^\top v = \vect(\nabla f_K)^\top J_K v = \|\Pi_{\mathcal{T}} \vect(\nabla f_K)\|^2$ and $\|v\|_G = \|\Pi_{\mathcal{T}} \vect(\nabla f_K)\|$, and the Cauchy--Schwarz inequality $\nabla F^\top v \leq \|\nabla F\|_{G^{-1}} \|v\|_G$ yields (ii). (iii) The Hessian of $\mathcal{B}$ and the Jacobian $J_K$ transform with $T$, whereas $\Pi_{\mathcal{T}}$ depends on $\theta$ only through $Z(\xi)$.
\end{proof}

Lemma~\ref{lem:metric} collects the properties for which $G$ is designed. By (i), \eqref{eq:hr} is well defined with $G$. Since $F$ decreases along \eqref{eq:hr} at the rate $\|\nabla F\|^2_{G^{-1}}$, (ii) ensures that $F$ decreases as long as $\Pi_{\mathcal{T}}\vect(\nabla f_K) \neq 0$, namely as long as $\nabla f_K Z^{-1} B_u$ violates the symmetry part of Proposition~\ref{prop:necessary}. By (iii), \eqref{eq:hr} with $G$ generates the same gains for every choice of $Z_p$ and $Z_1, \dots, Z_q$. For $v \in \mathcal{S}$, $\|v\|_G$ is the Frobenius norm of the induced gain change, so that the steepest descent restricted to $\mathcal{S}$ moves the gain along the orthogonal projection of $-\nabla f_K$ onto $\mathcal{T}$, as the policy gradient flow restricted to $\mathcal{T}$ does. A Euclidean completion $H + \Pi_{\mathcal{S}}$ with the orthogonal projector $\Pi_{\mathcal{S}}$ onto $\mathcal{S}$ is positive definite as well, but it depends on the coordinates of $W$, and the induced gain changes along $\mathcal{T}$ contain the factor $Z^{-1}$ of $\nabla_W F$ in addition to the factor $Z^{-1}$ of $K = W Z^{-1}$.

\subsection{Lifted Flow and Main Result}

The lifted Hessian--Riemannian gradient flow is
\begin{align}
    \dot\theta = -G(\theta)^{-1} \nabla F(\theta), \qquad \theta(0) = \theta_0 \in \Theta, \label{eq:flow}
\end{align}
with the gain $K(t) = W(t) Z(\xi(t))^{-1}$ and the storage function $S_t(x) = \tfrac{1}{2} x^\top Z(\xi(t))^{-1} x$. Let $\delta(\theta) \coloneqq \min\{ \lambda_{\min}(Z(\xi)), \lambda_{\min}(-\Nb(Z(\xi),W)) \}$ denote the margin of the certificate.

\begin{theorem}
    \label{thm:flow}
    Let Assumption~\ref{ass:sys} hold, $\theta_0 \in \Theta$, and $K_0 = K(0)$. Then the following statements hold.
    \begin{enumerate}
        \item[(i)] The flow \eqref{eq:flow} has a unique solution $\theta(t) \in \Theta$ for all $t \geq 0$.
        \item[(ii)] For all $t \geq 0$, $K(t) \in \mathcal{P}^+$, and \eqref{eq:closed} with $K = K(t)$ is strictly passive with the storage function $S_t$.
        \item[(iii)] $\tfrac{\mathrm{d}}{\mathrm{d}t} f_{K(t)} = -\|\nabla F(\theta(t))\|^2_{G(\theta(t))^{-1}} \leq 0$, with equality if and only if $K(t) = K^*$, and $K(t)$ remains in the compact set $\{ K \in \mathcal{K} \mid f_K \leq f_{K_0} \}$.
        \item[(iv)] For all $t \geq 0$,
        \begin{align}
            \delta(\theta(t)) \geq \delta(\theta_0) \, \exp\Bigl( -\sqrt{t \, (f_{K_0} - f_{K(t)})} \Bigr). \label{eq:margin}
        \end{align}
        \item[(v)] $\bar\theta \in \Theta$ is an equilibrium if and only if $\varphi(Z(\bar\xi),\bar W) = K^*$, and every accumulation point of $\theta(t)$ in $\Theta$ is an equilibrium. If $K^* \notin \mathcal{P}^+$, then $f_{K(t)}$ is strictly decreasing and $\theta(t)$ eventually leaves every compact subset of $\Theta$.
        \item[(vi)] $K(t)$ is invariant under every affine reparameterization $\theta = T \tilde\theta + c$ with invertible $T$ and does not depend on the choice of $U$ and $V$.
        \item[(vii)] If $K(t) \to \bar K$ and $Z(\xi(t)) \to \bar Z \succ 0$ as $t \to \infty$, then $\bar P = \bar Z^{-1}$ certifies $\bar K \in \mathcal{P}^0$, and $\nabla f_{\bar K} \bar P B_u$ is symmetric.
    \end{enumerate}
\end{theorem}

The proof is given in the Appendix. By (ii), the flow can be stopped at any time with a strictly passivating gain and its storage function. By (iii), the closed-loop systems along the flow are uniformly exponentially stable. By (v), the lifting introduces no spurious equilibria, and if $K^* \notin \mathcal{P}^+$, the lifted variable $\theta(t)$ approaches $\partial\Theta$ or diverges, while \eqref{eq:margin} excludes an exponential decay of the margin. By (vii), a limit with a nondegenerate storage matrix satisfies the symmetry part of the necessary condition of Proposition~\ref{prop:necessary}, which is the part that the completion \eqref{eq:G} adds for $m \geq 2$.

\begin{remark}
    For $m = 1$, \eqref{eq:flow} is the Hessian--Riemannian gradient flow of \cite{alvarez2004} with $h = \mathcal{B}$, and for $m \geq 2$, its metric is the barrier Hessian completed on the kernel $\mathcal{S}$. Neither the well-posedness result \cite[Theorem~4.1]{alvarez2004}, which requires $F$ to be defined on $\mathbb{R}^{n_\theta}$, nor the convergence results for convex objective functions apply, since $F$ is nonconvex and undefined where $Z$ is singular. Instead, Theorem~\ref{thm:flow}(i) follows from the coercivity of $f_K$, which keeps $K(t)$ in a compact set, and the logarithmic structure of $\mathcal{B}$, which yields \eqref{eq:margin}, whereas the convergence results are replaced by the weaker statements (iii), (v), and (vii).
\end{remark}

The right-hand side of \eqref{eq:flow} is the minimizer of $v \mapsto \nabla F(\theta)^\top v + \tfrac{1}{2}\|v\|^2_{G(\theta)}$, where the first term approximates $F(\theta + v) - F(\theta)$ to first order and the second term penalizes the length of $v$ in the metric $G$. Since $G$ grows unboundedly toward $\partial\Theta$, the flow slows down near the boundary, whereas in the interior it follows the descent of $F$. In contrast to the Euclidean gradient flow of the penalized objective $F + \tau\mathcal{B}$ with $\tau > 0$, whose equilibria are the stationary points of $F + \tau\mathcal{B}$ and whose trajectories depend on the coordinates of $\theta$, the equilibria of \eqref{eq:flow} are the stationary points of $F$, and its gain trajectories are invariant under affine changes of coordinates by (vi).

\subsection{Implementation}
\label{subsec:implementation}

\subsubsection*{Initialization}
Every solution $\theta_0$ of the strict LMIs $Z(\xi) \succ 0$ and $\Nb(Z(\xi),W) \prec 0$ of Lemma~\ref{lem:Theta}, which form a convex feasibility problem, is an admissible initial point, as is $(\xi_0, \vect(W))$ with $Z(\xi_0) = Z$ for every solution $(Z, W, \tilde N)$ of Proposition~\ref{prop:feas} with $\tilde N \succ 0$. If a strictly passivating gain $K_0$ is given without a certificate, the convex program
\begin{align}
    \max_{\xi, \delta} ~ & \delta - \epsilon \|Z(\xi)\|_F \notag \\
    \text{s.t.} ~ & Z(\xi) \succeq \delta I_n, ~ -\Nb(Z(\xi), K_0 Z(\xi)) \succeq \delta I_{n+r} \label{eq:init}
\end{align}
with $\epsilon > 0$ yields $\theta_0 = (\xi_0, \vect(K_0 Z(\xi_0))) \in \Theta$ whenever its optimal $\delta$ is positive, namely a certificate of $K_0$ with maximal margin up to the regularization. The regularization keeps the maximizer finite, since the feasible set of \eqref{eq:init} can be unbounded.

\subsubsection*{Evaluation Near the Boundary}
Near $\partial\Theta$, $\Gamma$ or $Z^{-1}$ in \eqref{eq:hessB} becomes unbounded, so that forming and inverting $G$ is inaccurate in floating-point arithmetic. In interior-point methods, where the same difficulty arises, standard factorizations of augmented systems often yield accurate search directions despite ill-conditioning \cite{Wright1997}, and we likewise compute $\dot\theta$ from an augmented system without forming $G$. With the eigendecompositions $-\Nb = O_N \diag(\lambda) O_N^\top$ and $Z = O_Z \diag(\nu) O_Z^\top$, let $\hat N_i \coloneqq O_N^\top (\partial \Nb / \partial \theta_i) O_N$ and $\hat Z_i \coloneqq O_Z^\top (\partial Z / \partial \theta_i) O_Z$ denote the derivatives in the eigenbases. Since $\Gamma = O_N \diag(\lambda)^{-1} O_N^\top$, the first term of \eqref{eq:hessB} is $\sum_{a,b} (\hat N_i)_{ab} (\hat N_k)_{ab} / (\lambda_a \lambda_b)$, where each pair $a \neq b$ appears twice by symmetry, and analogously for the second term with $\hat Z_i$ and $\nu$. Hence, $H = \tilde J^\top \diag(\omega)^{-1} \tilde J$, where the $i$-th column of $\tilde J$ contains the entries $(\hat N_i)_{ab}$ and $(\hat Z_i)_{ab}$ with $a \leq b$, and $\omega$ contains the corresponding weights $\omega_{ab} = \lambda_a \lambda_b$ for $a = b$ and $\omega_{ab} = \lambda_a\lambda_b/2$ for $a < b$, with $\nu$ in place of $\lambda$ for the entries of $\hat Z_i$. With an orthonormal basis $O_{\mathcal{T}}$ of $\vect(\mathcal{T})$ and $J_{\mathcal{T}} = O_{\mathcal{T}}^\top J_K$, it follows that $G = \tilde J^\top \diag(\omega)^{-1} \tilde J + J_{\mathcal{T}}^\top J_{\mathcal{T}}$, where $\tilde J$ is bounded and some weights in $\omega$ tend to zero at $\partial\Theta$, and, at each evaluation of the right-hand side of \eqref{eq:flow}, $\dot\theta$ is the first block of the solution of
\begin{align}
    \begin{bmatrix} 0 & \tilde J^\top & J_{\mathcal{T}}^\top \\ \tilde J & -\diag(\omega) & 0 \\ J_{\mathcal{T}} & 0 & -I \end{bmatrix} \begin{bmatrix} \dot\theta \\ \eta \\ \zeta \end{bmatrix} = \begin{bmatrix} -\nabla F(\theta) \\ 0 \\ 0 \end{bmatrix}, \label{eq:augmented}
\end{align}
whose auxiliary variables $\eta = \diag(\omega)^{-1}\tilde J\dot\theta$ and $\zeta = J_{\mathcal{T}}\dot\theta$ recover $G\dot\theta = -\nabla F(\theta)$ upon elimination, and whose entries remain bounded as $\Nb$ becomes singular while $Z$ remains nonsingular.

\subsubsection*{Computational Effort}
Each evaluation requires two Lyapunov equations of size $n$, two symmetric eigendecompositions of sizes $n + r$ and $n$, an orthonormal basis of $\mathcal{T}$ for $m \geq 2$, and a linear solve \eqref{eq:augmented} of size $n_\theta + (n+r)(n+r+1)/2 + n(n+1)/2 + m(m-1)/2$ with $n_\theta \leq n(n+1)/2 + mn$, whose dense factorization dominates the cost.

\section{Numerical Results}
\label{sec:numerics}

Section~\ref{subsec:compare} compares the flow with the indirect approach, Section~\ref{subsec:ft} considers nonconvex passivity regions, and Section~\ref{subsec:multi} considers multi-input systems. In all examples, the flow is integrated with the MATLAB solver \texttt{ode15s}. The examples of Sections~\ref{subsec:compare} and~\ref{subsec:ft} have $n = 2$, $m = p = 1$, and $Q = I_2$, so that $G = H$, the flow is initialized by \eqref{eq:init} with $\epsilon = 10^{-3}$, and the reference $\Kopt$ is a minimizer of $f_K$ over $\mathcal{K} \cap \mathcal{P}^0$, computed with the MATLAB function \texttt{fmincon} from an analytic description of $\mathcal{P}^0$, so that $f_{\Kopt}$ is the infimum of $f_K$ over $\mathcal{P}^+$ by Lemma~\ref{lem:closure}.

\subsection{Comparison with the Indirect Approach}
\label{subsec:compare}

We consider the open-loop stable system with $D = 0$, $R = 2$, and
\begin{align}
    A = \begin{tbmatrix} -2 & -1 \\ -1 & -3 \end{tbmatrix}, ~ B_u = \begin{tbmatrix} 1 \\ 2 \end{tbmatrix}, ~ B_d = \begin{tbmatrix} 2 \\ 1 \end{tbmatrix}, ~ C = \begin{tbmatrix} 0 & 1 \end{tbmatrix}, \label{eq:sys1}
\end{align}
and the open-loop unstable system with $D = 0$, $R = 1$, and
\begin{align}
    A = \begin{tbmatrix} 0 & 1 \\ -3 & 1 \end{tbmatrix}, ~ B_u = \begin{tbmatrix} 0 \\ 1 \end{tbmatrix}, ~ B_d = \begin{tbmatrix} -1 \\ 2 \end{tbmatrix}, ~ C = \begin{tbmatrix} 1 & 1 \end{tbmatrix}, \label{eq:sys2}
\end{align}
whose passivity regions are polygons given in \cite{giessler2026passive}. For both systems, $q = 1$ and $n_\theta = 3$, the flow is evaluated by solving with $H$ directly, and it is integrated until $\|\nabla F\|_{H^{-1}} < 10^{-8}$ or $t = 10^7$. Table~\ref{tab:results} compares the end points of the flow from three initial gains, which agree to the digits shown, with $K^*$ and the gain $\hat K$ of the indirect approach with cubes of edge length $0.1$ \cite{giessler2026passive}.

\begin{table}[t]
    \centering
    \caption{Gains and relative cost increase over $f_{K^*}$.}
    \label{tab:results}
    \setlength{\tabcolsep}{3pt}
    \footnotesize
    \begin{tabular}{lcccc}
        \toprule
        & \multicolumn{2}{c}{System \eqref{eq:sys1}} & \multicolumn{2}{c}{System \eqref{eq:sys2}} \\
        \cmidrule(lr){2-3} \cmidrule(lr){4-5}
        & $K$ & increase & $K$ & increase \\
        \midrule
        $K^*$ & $[0.0484 ~~ 0.1430]$ & $0$ & $[0.1623 ~~ 2.5246]$ & $0$ \\
        flow \eqref{eq:flow} & $[0.0000 ~~ 0.1172]$ & $0.24\,\%$ & $[0.3882 ~~ 3.6941]$ & $8.53\,\%$ \\
        $\hat K$ \cite{giessler2026passive} & $[-0.1000 ~~ 0.0694]$ & $2.41\,\%$ & $[0.4467 ~~ 3.8234]$ & $9.77\,\%$ \\
        \bottomrule
    \end{tabular}
\end{table}

\begin{figure*}[t]
    \centering
    \includegraphics{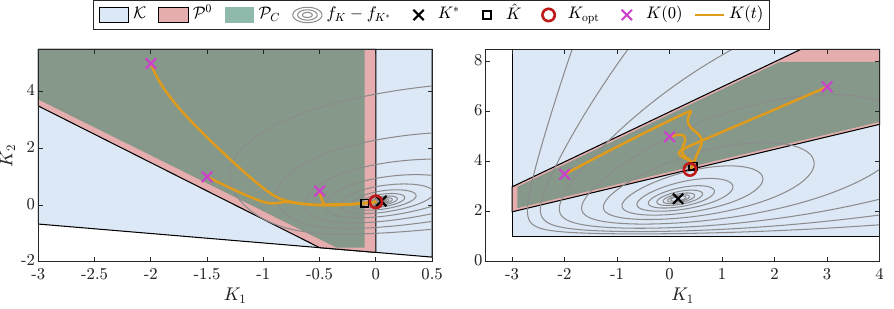}
    \caption{System~\eqref{eq:sys1} (left) and system~\eqref{eq:sys2} (right). Stability region $\mathcal{K}$, passivity region $\mathcal{P}^0$, inner approximation $\mathcal{P}_C$ of the indirect approach, level sets of $f_K - f_{K^*}$, the gains $K^*$, $\hat K$, and $\Kopt$, and trajectories $K(t)$ of the lifted flow from three initial gains $K(0)$.}
    \label{fig:sys}
\end{figure*}

In Fig.~\ref{fig:sys}, the flow approaches $\Kopt$ from all initial gains, up to a distance of $2.3 \cdot 10^{-6}$ for \eqref{eq:sys1} and $1.1 \cdot 10^{-7}$ for \eqref{eq:sys2}, and the cost decreases monotonically. Since $\Kopt$ lies on $\partial\mathcal{P}^+$ outside the inner approximation $\mathcal{P}_C$, to which the indirect approach restricts $\hat K$, the flow attains a lower cost than $\hat K$. The margin tends to zero, and for \eqref{eq:sys2}, the cost gap and the margin decay like $1/t$, with fitted log-log slopes of $-1.00$ over $t \in [10^4, 10^7]$.

\subsection{Nonconvex Passivity Regions}
\label{subsec:ft}

The six systems of Table~\ref{tab:ft} have $D = 0.1$ and $R = 1$, so that $r = 1$, $V$ is empty, and $n_\theta = 5$. With the closed-loop transfer function $D + (b_1 s + b_0)/(s^2 + a_1 s + a_0)$, where $a_1 = -\tr(A_K)$, $a_0 = \det(A_K)$, $b_1 = CB_d$, and $b_0 = CA_KB_d + a_1 b_1$, the positive realness conditions \cite[Definition~6.4]{khalil_nonlinear_2002} yield
\begin{align*}
    \mathcal{P}^0 \cap \mathcal{K} = \bigl\{ K \in \mathcal{K} \bigm| D a_0 + b_0 \geq 0, ~ \beta_1 + 2\sqrt{D \beta_0} \geq 0 \bigr\}
\end{align*}
with $\beta_1 = D(a_1^2 - 2a_0) + a_1 b_1 - b_0$ and $\beta_0 = a_0 (D a_0 + b_0)$, which is nonconvex for all six systems. The flow is evaluated by \eqref{eq:augmented}, started from four initial gains per system spread over $\operatorname{int}(\mathcal{P}^0)$, and integrated until $\|\nabla F\|_{G^{-1}} < 10^{-10}$, $\delta < 10^{-10}$, or $t = 10^{10}$.

\begin{table}[t]
    \centering
    \caption{Systems with $D = 0.1$, minimizer $\Kopt$ with relative cost increase over $f_{K^*}$, and number of initial gains from which the flow converges to $\Kopt$.}
    \label{tab:ft}
    \setlength{\tabcolsep}{2pt}
    \renewcommand{\arraystretch}{1.5}
    \footnotesize
    \begin{tabular}{lccccccc}
        \toprule
        & $A$ & $B_u^\top$ & $B_d^\top$ & $C$ & $\Kopt$ & incr. & conv. \\
        \midrule
        I   & $\bigl[\begin{smallmatrix} -3 & -3 \\ 2 & 2 \end{smallmatrix}\bigr]$ & $[-1 ~ 2]$ & $[2 ~ {-2}]$ & $[2 ~ 2]$ & $[1.306 ~ 1.830]$ & $26.9\,\%$ & 2/4 \\
        II  & $\bigl[\begin{smallmatrix} -3 & 3 \\ 1 & -3 \end{smallmatrix}\bigr]$ & $[1 ~ {-1}]$ & $[0 ~ 2]$ & $[2 ~ 0]$ & $[3.400 ~ {-2.774}]$ & $178\,\%$ & 4/4 \\
        III & $\bigl[\begin{smallmatrix} 0 & 3 \\ 2 & 2 \end{smallmatrix}\bigr]$ & $[-1 ~ {-2}]$ & $[2 ~ 1]$ & $[1 ~ {-2}]$ & $[-2.043 ~ {-5.367}]$ & $17.6\,\%$ & 4/4 \\
        IV  & $\bigl[\begin{smallmatrix} -1 & 0 \\ 2 & 0 \end{smallmatrix}\bigr]$ & $[2 ~ 1]$ & $[0 ~ {-2}]$ & $[2 ~ 0]$ & $[2.492 ~ 1.000]$ & $31.9\,\%$ & 4/4 \\
        V   & $\bigl[\begin{smallmatrix} -2 & -2 \\ 2 & -2 \end{smallmatrix}\bigr]$ & $[-1 ~ {-1}]$ & $[0 ~ {-2}]$ & $[2 ~ 0]$ & $[-2.854 ~ {-0.562}]$ & $195\,\%$ & 4/4 \\
        VI  & $\bigl[\begin{smallmatrix} -3 & 0 \\ 0 & 0 \end{smallmatrix}\bigr]$ & $[-2 ~ {-1}]$ & $[0 ~ {-2}]$ & $[2 ~ 0]$ & $[-2.127 ~ {-1.364}]$ & $26.7\,\%$ & 4/4 \\
        \bottomrule
    \end{tabular}
\end{table}

\begin{figure*}[t]
    \centering
    \includegraphics{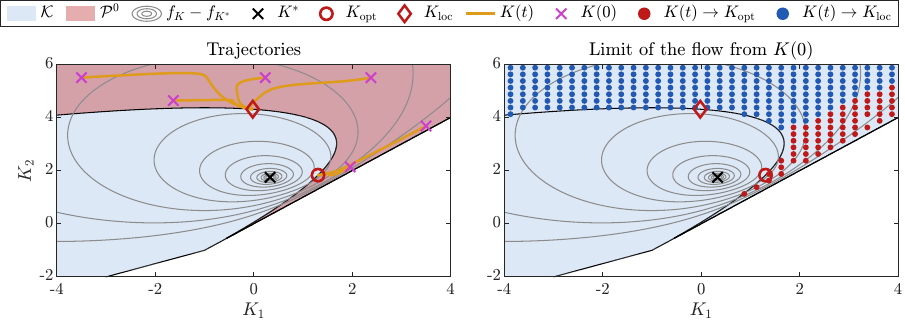}
    \caption{System~I of Table~\ref{tab:ft}. Left, stability region $\mathcal{K}$, passivity region $\mathcal{P}^0$, level sets of $f_K - f_{K^*}$, $K^*$, the global minimizer $\Kopt$ and the local minimizer $\Kloc$ of $f_K$ over $\mathcal{K} \cap \mathcal{P}^0$, and trajectories $K(t)$ of the lifted flow from six initial gains $K(0)$. Right, initial gains on a grid with spacing $0.25$, marked by the minimizer to which the flow converges.}
    \label{fig:nc}
\end{figure*}

In Systems~II to~VI, all trajectories converge to $\Kopt$ with a relative cost gap below $4.0 \cdot 10^{-8}$. In System~I, $f_K$ has a second local minimizer $\Kloc = [-0.020 ~~ 4.311]$ on $\partial\mathcal{P}^+$ with a cost increase of $36.8\,\%$, and two of the four trajectories converge to it. Fig.~\ref{fig:nc} shows System~I in detail. On the left, two of six trajectories converge to $\Kopt$ and four to $\Kloc$, and the cost gap decays like $1/t$ toward both limits, so that the convergence behavior does not indicate whether the limit is the global minimizer. On the right, of the $284$ initial gains on a grid with spacing $0.25$ in $\operatorname{int}(\mathcal{P}^0)$, with stopping tolerances relaxed to $10^{-6}$, $58$ lead to $\Kopt$ and $224$ to $\Kloc$, whereas the initialization finds no certificate for two gains close to $\partial\mathcal{P}^0$.

\subsection{Multi-Input Systems}
\label{subsec:multi}

The first system has $n = 3$, $m = 2$, $p = 1$, $D = 0$, $B_d^\top = [1 ~~ {-2} ~~ {-2}]$, $C = [6 ~~ 0 ~~ {-7}]$, and
\begin{align}
    A = \begin{tbmatrix} -2 & 2 & -3 \\ -3 & -1 & 0 \\ 0 & 0 & 3 \end{tbmatrix}, ~ B_u = \begin{tbmatrix} 2 & 2 \\ 1 & 2 \\ -2 & -2 \end{tbmatrix}, \label{eq:sys3}
\end{align}
so that $n_\theta = 9$ and $\dim \mathcal{S} = 1$. The second system has $n = 4$, $m = 3$, $p = 2$, the singular $D + D^\top \neq 0$ with $D = \diag(0.2, 0)$, and
\begin{align}
    A &= \begin{tbmatrix} 0 & 3 & -3 & 1 \\ 1 & 3 & 0 & -2 \\ -1 & -3 & -3 & 1 \\ 0 & -2 & 0 & 1 \end{tbmatrix}, \quad B_u = \begin{tbmatrix} -2 & -1 & 0 \\ 0 & 1 & 1 \\ -1 & 0 & -1 \\ 0 & -2 & -1 \end{tbmatrix}, \notag \\
    B_d^\top &= \begin{tbmatrix} -1 & 2 & 2 & 1 \\ 2 & 2 & -1 & -2 \end{tbmatrix}, \quad C = \begin{tbmatrix} 1 & 3 & 7 & 6 \\ 1 & 6 & -5 & -12 \end{tbmatrix}, \label{eq:sys4}
\end{align}
so that $r = 1$, $\Nb$ has size $5$, $n_\theta = 18$, and $\dim \mathcal{S} = 3$. Both systems use $Q = I_n$ and $R = I_m$, and both are open-loop unstable. The five initial points per system solve the LMIs of Lemma~\ref{lem:Theta} with $\|Z\|_F, \|W\|_F \leq 20$, namely the one with the largest margin $\delta_{\max}$ and the maximizers of four linear functions of $W$ with standard normally distributed coefficients subject to a margin of at least $0.3\,\delta_{\max}$, so that the initial gains are spread over $\mathcal{P}^+$. From each initial point, three variants of the flow are integrated, namely with the gain completion \eqref{eq:G}, with the Euclidean completion $H + \Pi_{\mathcal{S}}$ in place of $G$, and with no completion, $\dot\theta = -H^\dagger \nabla F$, each evaluated by an augmented system as in \eqref{eq:augmented}, until the criteria of Section~\ref{subsec:ft} are met or $300\,$s of computation time have elapsed. The reference $\Kopt$ is the best result of \texttt{fmincon} on the closure of $\Theta$ from $20$ initial points and from the end points of the flow. It increases the cost over $f_{K^*}$ by $9.36\,\%$ for \eqref{eq:sys3} and by $30.30\,\%$ for \eqref{eq:sys4}.

\begin{figure*}[t]
    \centering
    \includegraphics{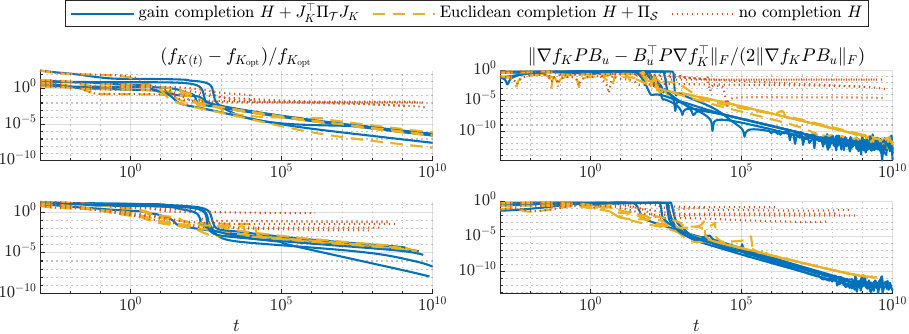}
    \caption{Multi-input systems~\eqref{eq:sys3} (top) and~\eqref{eq:sys4} (bottom). Relative cost gap with respect to the reference (left) and relative skew-symmetric part of $\nabla f_K P B_u$ with $P = Z^{-1}$ (right) along the flow from five initial points, for the gain completion \eqref{eq:G}, the Euclidean completion, and no completion.}
    \label{fig:multi}
\end{figure*}

Fig.~\ref{fig:multi} shows the relative cost gap (left) and the relative skew-symmetric part of $\nabla f_K Z^{-1} B_u$ (right) along all trajectories. With the gain completion, the flow approaches $\Kopt$ from all initial points. The relative cost gap at the end points is at most $5.2 \cdot 10^{-7}$ for \eqref{eq:sys3} and at most $4.7 \cdot 10^{-5}$ for \eqref{eq:sys4}, where three of the five trajectories reach the time limit. The cost decreases monotonically, the margin remains positive, and at the end points, the relative skew-symmetric part of $\nabla f_K Z^{-1} B_u$ is at most $2.7 \cdot 10^{-12}$ and the relative smallest eigenvalue of its symmetric part is above $-2.0 \cdot 10^{-6}$, in accordance with Proposition~\ref{prop:necessary} and Theorem~\ref{thm:flow}(vii). The Euclidean completion behaves similarly for \eqref{eq:sys3}, whereas for \eqref{eq:sys4}, it reduces the cost gap and the skew-symmetric part more slowly, and all five trajectories reach the time limit with a relative cost gap of up to $7.4 \cdot 10^{-5}$. Without completion, the flow stalls between $0.2\,\%$ and $1.6\,\%$ above $f_{\Kopt}$ for \eqref{eq:sys3} and between $0.6\,\%$ and $90\,\%$ above $f_{\Kopt}$ for \eqref{eq:sys4}, with a skew-symmetric part of up to $0.17$, since the lifted variable cannot move along $\mathcal{S}$.

\section{Discussion and Limitations}
\label{sec:discussion}

In contrast to the indirect approach, which certifies a polytope of gains by a common storage function and is therefore conservative and restricted to few gain entries, the lifted flow certifies each gain individually, which removes the conservatism and the exponential growth with $mn$. Since this certificate is strict, it also covers a neighborhood of the current gain that shrinks only as the margin $\delta$ tends to zero, so that stopping the flow at a prescribed margin trades cost for robustness against gain perturbations.

The flow is a local method. Theorem~\ref{thm:flow} guarantees strict passivity, descent, and global existence, but not convergence to a minimizer of \eqref{eq:problem}, and System~I shows that the passivity-constrained LQR problem can have local minimizers, in contrast to the unconstrained LQR problem and to the problems with convex parameterizations of the cost \cite{Fazel2021,zheng_benign_2026}. Solvers for nonlinear programs, such as \texttt{fmincon} used for the references $\Kopt$, are local as well, and their iterates satisfy the constraints in general only at convergence, whereas every gain along the flow is strictly passivating with an explicit storage function, so that the flow can be stopped at any time. A characterization of the limit points on $\partial\mathcal{P}^+$, for example as Karush--Kuhn--Tucker points of \eqref{eq:problem}, remains open.

The approach requires $\mathcal{P}^+ \neq \emptyset$ in Assumption~\ref{ass:sys}, since the barrier is defined only if some lifted gain satisfies $\Nb \prec 0$, and it thus excludes systems that state feedback renders passive but not strictly passive, such as $A = \begin{bsmallmatrix} 0 & 1 \\ a_{21} & a_{22} \end{bsmallmatrix}$, $B_u = e_2$, $B_d = C^\top = e_1$, and $D = 0$, for which $PB_d = C^\top$ implies $(A_K^\top P + PA_K)_{11} = 0$ for all $K$.

Since the metric diverges on $\partial\Theta$, the convergence to boundary minimizers is slow, and the barrier Hessian becomes ill-conditioned, which the augmented system \eqref{eq:augmented} handles down to margins of $10^{-10}$ in the examples. For \eqref{eq:sys4}, three of the five trajectories reach the limit of $300\,$s of computation time with a relative cost gap of up to $4.7 \cdot 10^{-5}$. The behavior for systems of higher order remains to be evaluated.

\section{Conclusion and Outlook}
\label{sec:conclusion}

We proposed a direct approach to LQR synthesis under passivity constraints, based on an explicit description of the convex lifted set of strictly passivating gains for every admissible feedthrough matrix and on a Hessian--Riemannian gradient flow with a logarithmic barrier as metric. For multi-input systems, the barrier Hessian is singular along certificate-preserving gain changes, and its completion by the gain metric leads to limits that satisfy the symmetry part of a new necessary condition. The flow exists for all times, generates strictly passivating gains with explicit storage functions, and decreases the cost monotonically. In the examples, it attained a lower cost than the indirect approach, approached the reference minimizers of two multi-input systems, one of them with singular symmetric feedthrough, and on a nonconvex passivity region, it converged to a local minimizer from part of the initial gains. Future work includes discretizations with convergence guarantees, a characterization of the limit points, and extensions to other dissipativity constraints with quadratic supply rates and to data-driven formulations.

\appendix
\section*{Proof of Theorem~\ref{thm:flow}}

\subsubsection*{Local existence, (ii), and (iii)}
By Lemma~\ref{lem:gradF}, \eqref{eq:G}, and $G \succ 0$ from Lemma~\ref{lem:metric}(i), the right-hand side of \eqref{eq:flow} is real analytic and thus locally Lipschitz on $\Theta$ \cite[Lemma~3.2]{khalil_nonlinear_2002}, so that a unique solution exists on a maximal interval $[0, T_{\max})$ with $T_{\max} \in (0, \infty]$ \cite[Theorem~3.1 and p.~93]{khalil_nonlinear_2002}. On this interval, $(Z(\xi(t)), W(t)) \in \mathcal{L}^+$, which shows (ii) by Proposition~\ref{prop:feas} and Lemma~\ref{lem:KYP}. The chain rule gives $\tfrac{\mathrm{d}}{\mathrm{d}t} F(\theta) = -\nabla F^\top G^{-1} \nabla F$, which vanishes if and only if $K(t) = K^*$ by Lemma~\ref{lem:gradF}, and $K(t)$ remains in the sublevel set, which is compact by Section~\ref{subsec:lqr}.

\subsubsection*{Statement (iv)}
Let $\Upsilon(t)$ denote $Z(\xi(t))$ or $-\Nb(Z(\xi(t)), W(t))$, both positive definite, and let $v$ be a fixed unit vector. The function $g = v^\top \Upsilon v$ satisfies $\dot g / g = \hat v^\top \Upsilon^{-1/2} \dot\Upsilon \Upsilon^{-1/2} \hat v$ with the unit vector $\hat v = \Upsilon^{1/2} v / \|\Upsilon^{1/2} v\|$, so that
\begin{align*}
    |\dot g / g| \leq \tr(\Upsilon^{-1} \dot\Upsilon \Upsilon^{-1} \dot\Upsilon)^{1/2} \leq \|\dot\theta\|_{H} \leq \|\dot\theta\|_{G},
\end{align*}
since $\tr(\Upsilon^{-1}\dot\Upsilon \Upsilon^{-1} \dot\Upsilon)$ is one of the two nonnegative terms of $\dot\theta^\top H \dot\theta$ by \eqref{eq:hessB}, and $G \succeq H$ by Lemma~\ref{lem:metric}. With $\|\dot\theta\|_G^2 = \nabla F^\top G^{-1} \nabla F = -\tfrac{\mathrm{d}}{\mathrm{d}t} F(\theta)$ and the Cauchy--Schwarz inequality, $\int_0^t \|\dot\theta(\tau)\|_G \, \mathrm{d}\tau \leq \sqrt{t \, (f_{K_0} - f_{K(t)})} \eqqcolon \ell(t)$. Hence $g(0) e^{-\ell(t)} \leq g(t) \leq g(0) e^{\ell(t)}$ for every unit vector $v$, which yields $\lambda_{\min}(\Upsilon(t)) \geq \lambda_{\min}(\Upsilon(0)) e^{-\ell(t)}$ for both choices of $\Upsilon$, and thus \eqref{eq:margin}, and $\lambda_{\max}(\Upsilon(t)) \leq \lambda_{\max}(\Upsilon(0)) e^{\ell(t)}$.

\subsubsection*{Statement (i)}
A solution can cease to exist at a finite time only by diverging or by approaching the boundary of the domain of the right-hand side, namely, if $T_{\max} < \infty$, then $\theta(t)$ leaves every compact subset of $\Theta$ \cite[p.~93]{khalil_nonlinear_2002}. Suppose $T_{\max} < \infty$. Then $\ell(t) \leq \bar\ell \coloneqq (T_{\max} f_{K_0})^{1/2}$ on $[0, T_{\max})$, since $f_K \geq 0$, and the bounds in the proof of (iv) yield $\delta(\theta(t)) \geq \underline\delta \coloneqq \delta(\theta_0) e^{-\bar\ell} > 0$, so that $\theta(t)$ does not approach $\partial\Theta$, and $\lambda_{\max}(Z(\xi(t))) \leq \lambda_{\max}(Z(\xi_0)) e^{\bar\ell}$. Since $\xi \mapsto Z(\xi)$ is affine and injective, $\xi(t)$ is thus bounded, and since $K(t)$ is bounded by (iii), so is $W(t) = K(t) Z(\xi(t))$, so that $\theta(t)$ does not diverge. Hence, $\theta(t)$ remains in the set $\{\theta \mid \|\theta\| \leq c,~ \delta(\theta) \geq \underline\delta\}$ for some $c > 0$, which is a compact subset of $\Theta$, since $\delta$ is continuous. This contradicts $T_{\max} < \infty$, so that $T_{\max} = \infty$.

\subsubsection*{Statement (v)}
Since $G \succ 0$, the equilibria are the zeros of $\nabla F$, which correspond to $K^*$ by Lemma~\ref{lem:gradF}. Let $\theta(t_k) \to \bar\theta \in \Theta$ with $t_k \to \infty$ and $\nabla F(\bar\theta) \neq 0$. By continuity, there exist $\rho, c_1, c_2 > 0$ such that the closed ball with center $\bar\theta$ and radius $2\rho$ lies in $\Theta$ and $\|\nabla F\|^2_{G^{-1}} \geq c_1$ and $\|G^{-1} \nabla F\| \leq c_2$ hold on it. For every $k$ with $\|\theta(t_k) - \bar\theta\| \leq \rho$, the trajectory remains in this ball on $[t_k, t_k + \rho/c_2]$, so that $F$ decreases by at least $c_1\rho/c_2$ on this interval. Passing to a subsequence with $t_{k+1} \geq t_k + \rho/c_2$ yields infinitely many disjoint such intervals, which contradicts $F \geq 0$. If $K^* \notin \mathcal{P}^+$, no equilibrium and hence no accumulation point exists in $\Theta$. Then $\theta(t)$ eventually leaves every compact set $\mathcal{X} \subseteq \Theta$, since a sequence $\theta(t_k) \in \mathcal{X}$ with $t_k \to \infty$ would have an accumulation point in $\mathcal{X}$, and $f_{K(t)}$ is strictly decreasing by (iii), since $K(t) \in \mathcal{P}^+$ implies $K(t) \neq K^*$.

\subsubsection*{Statement (vi)}
In the coordinates $\tilde\theta$, the gradient is $T^\top \nabla F$ and the metric is $T^\top G T$ by Lemma~\ref{lem:metric}(iii). Hence $\dot{\tilde\theta} = -(T^\top G T)^{-1} T^\top \nabla F = -T^{-1} G^{-1} \nabla F$, namely $T \dot{\tilde\theta} = \dot\theta$, so that both coordinates generate the same $(Z(\xi(t)), W(t))$. Another choice of $U$ and $V$ leaves the affine subspace, $\log\det(-\Nb)$, $J_K$, and $\Pi_{\mathcal{T}}$ unchanged.

\subsubsection*{Statement (vii)}
By (iii) and Lemma~\ref{lem:metric}(ii), $\int_0^\infty \|\Pi_{\mathcal{T}} \vect(\nabla f_{K(t)})\|^2 \mathrm{d}t \leq \int_0^\infty \|\nabla F\|^2_{G^{-1}} \mathrm{d}t \leq f_{K_0}$. Since $\bar K \in \mathcal{K}$ by (iii) and $\Pi_{\mathcal{T}}$ depends continuously on $Z \succ 0$, the integrand converges to $\|\Pi_{\mathcal{T}(\bar P)} \vect(\nabla f_{\bar K})\|^2$, which therefore vanishes, namely $\tr(\Sigma^\top \nabla f_{\bar K} \bar P B_u) = 0$ for all skew-symmetric $\Sigma$, so that $\nabla f_{\bar K} \bar P B_u$ is symmetric. Moreover, $(\bar Z, \bar K \bar Z)$ satisfies the conditions of \eqref{eq:L0} as a limit of points of $\mathcal{L}^+$ with $\bar Z \succ 0$, so that $\bar P$ certifies $\bar K \in \mathcal{P}^0$ by Theorem~\ref{thm:lifted} and Proposition~\ref{prop:feas}.

\printbibliography

\end{document}